\documentclass{amsart}

\input{epsf.tex}
\usepackage{amsfonts,amssymb,verbatim,amsmath,amsthm,latexsym,textcomp,amscd}
\usepackage{latexsym,amsfonts,amssymb,epsfig,verbatim}
\usepackage{amsmath,amsthm,amssymb,latexsym,graphics,textcomp}
\usepackage{graphicx}
\usepackage{color}
\usepackage{url}

\begin{document}

\newtheorem{theorem}{Theorem}[section]
\newtheorem{prop}[theorem]{Proposition}
\newtheorem{lemma}[theorem]{Lemma}
\newtheorem{cor}[theorem]{Corollary}
\newtheorem{defn}[theorem]{Definition}
\newtheorem{conj}[theorem]{Conjecture}
\newtheorem{claim}[theorem]{Claim}

\newcommand{\boundary}{\partial}
\newcommand{\bbC}{{\mathbb C}}
\newcommand{\bbD}{{\mathbb D}}
\newcommand{\bbH}{{\mathbb H}}
\newcommand{\bbZ}{{\mathbb Z}}
\newcommand{\bbN}{{\mathbb N}}
\newcommand{\bbQ}{{\mathbb Q}}
\newcommand{\bbR}{{\mathbb R}}
\newcommand{\proj}{{\mathbb P}}
\newcommand{\lhp}{{\mathbb L}}
\newcommand{\tube}{{\mathbb T}}
\newcommand{\cusp}{{\mathbb P}}
\newcommand\AAA{{\mathcal A}}
\newcommand\BB{{\mathcal B}}
\newcommand\CC{{\mathcal C}}
\newcommand\DD{{\mathcal D}}
\newcommand\EE{{\mathcal E}}
\newcommand\FF{{\mathcal F}}
\newcommand\GG{{\mathcal G}}
\newcommand\HH{{\mathcal H}}
\newcommand\II{{\mathcal I}}
\newcommand\JJ{{\mathcal J}}
\newcommand\KK{{\mathcal K}}
\newcommand\LL{{\mathcal L}}
\newcommand\MM{{\mathcal M}}
\newcommand\NN{{\mathcal N}}
\newcommand\OO{{\mathcal O}}
\newcommand\PP{{\mathcal P}}
\newcommand\QQ{{\mathcal Q}}
\newcommand\RR{{\mathcal R}}
\newcommand\SSS{{\mathcal S}}
\newcommand\TT{{\mathcal T}}
\newcommand\UU{{\mathcal U}}
\newcommand\VV{{\mathcal V}}
\newcommand\WW{{\mathcal W}}
\newcommand\XX{{\mathcal X}}
\newcommand\YY{{\mathcal Y}}
\newcommand\ZZ{{\mathcal Z}}
\newcommand\CH{{\CC\HH}}
\newcommand\TC{{\TT\CC}}
\newcommand\EXH{{ \EE (X, \HH )}}
\newcommand\GXH{{ \GG (X, \HH )}}
\newcommand\GYH{{ \GG (Y, \HH )}}
\newcommand\PEX{{\PP\EE  (X, \HH , \GG , \LL )}}
\newcommand\MF{{\MM\FF}}
\newcommand\PMF{{\PP\kern-2pt\MM\FF}}
\newcommand\ML{{\MM\LL}}
\newcommand\PML{{\PP\kern-2pt\MM\LL}}
\newcommand\GL{{\GG\LL}}
\newcommand\Pol{{\mathcal P}}
\newcommand\half{{\textstyle{\frac12}}}
\newcommand\Half{{\frac12}}
\newcommand\Mod{\operatorname{Mod}}
\newcommand\Area{\operatorname{Area}}
\newcommand\ep{\epsilon}
\newcommand\hhat{\widehat}
\newcommand\Proj{{\mathbf P}}
\newcommand\U{{\mathbf U}}
 \newcommand\Hyp{{\mathbf H}}
\newcommand\D{{\mathbf D}}
\newcommand\Z{{\mathbb Z}}
\newcommand\R{{\mathbb R}}
\newcommand\Q{{\mathbb Q}}
\newcommand\E{{\mathbb E}}
\newcommand\til{\widetilde}
\newcommand\length{\operatorname{length}}
\newcommand\tr{\operatorname{tr}}
\newcommand\gesim{\succ}
\newcommand\lesim{\prec}
\newcommand\simle{\lesim}
\newcommand\simge{\gesim}
\newcommand{\simmult}{\asymp}
\newcommand{\simadd}{\mathrel{\overset{\text{\tiny $+$}}{\sim}}}
\newcommand{\ssm}{\setminus}
\newcommand{\diam}{\operatorname{diam}}
\newcommand{\pair}[1]{\langle #1\rangle}
\newcommand{\T}{{\mathbf T}}
\newcommand{\inj}{\operatorname{inj}}
\newcommand{\pleat}{\operatorname{\mathbf{pleat}}}
\newcommand{\short}{\operatorname{\mathbf{short}}}
\newcommand{\vertices}{\operatorname{vert}}
\newcommand{\collar}{\operatorname{\mathbf{collar}}}
\newcommand{\bcollar}{\operatorname{\overline{\mathbf{collar}}}}
\newcommand{\I}{{\mathbf I}}
\newcommand{\tprec}{\prec_t}
\newcommand{\fprec}{\prec_f}
\newcommand{\bprec}{\prec_b}
\newcommand{\pprec}{\prec_p}
\newcommand{\ppreceq}{\preceq_p}
\newcommand{\sprec}{\prec_s}
\newcommand{\cpreceq}{\preceq_c}
\newcommand{\cprec}{\prec_c}
\newcommand{\topprec}{\prec_{\rm top}}
\newcommand{\Topprec}{\prec_{\rm TOP}}
\newcommand{\fsub}{\mathrel{\scriptstyle\searrow}}
\newcommand{\bsub}{\mathrel{\scriptstyle\swarrow}}
\newcommand{\fsubd}{\mathrel{{\scriptstyle\searrow}\kern-1ex^d\kern0.5ex}}
\newcommand{\bsubd}{\mathrel{{\scriptstyle\swarrow}\kern-1.6ex^d\kern0.8ex}}
\newcommand{\fsubeq}{\mathrel{\raise-.7ex\hbox{$\overset{\searrow}{=}$}}}
\newcommand{\bsubeq}{\mathrel{\raise-.7ex\hbox{$\overset{\swarrow}{=}$}}}
\newcommand{\tw}{\operatorname{tw}}
\newcommand{\base}{\operatorname{base}}
\newcommand{\trans}{\operatorname{trans}}
\newcommand{\rest}{|_}
\newcommand{\bbar}{\overline}
\newcommand{\UML}{\operatorname{\UU\MM\LL}}
\newcommand{\EL}{\mathcal{EL}}
\newcommand{\tsum}{\sideset{}{'}\sum}
\newcommand{\tsh}[1]{\left\{\kern-.9ex\left\{#1\right\}\kern-.9ex\right\}}
\newcommand{\Tsh}[2]{\tsh{#2}_{#1}}
\newcommand{\qeq}{\mathrel{\approx}}
\newcommand{\Qeq}[1]{\mathrel{\approx_{#1}}}
\newcommand{\qle}{\lesssim}
\newcommand{\Qle}[1]{\mathrel{\lesssim_{#1}}}
\newcommand{\simp}{\operatorname{simp}}
\newcommand{\vsucc}{\operatorname{succ}}
\newcommand{\vpred}{\operatorname{pred}}
\newcommand\fhalf[1]{\overrightarrow {#1}}
\newcommand\bhalf[1]{\overleftarrow {#1}}
\newcommand\sleft{_{\text{left}}}
\newcommand\sright{_{\text{right}}}
\newcommand\sbtop{_{\text{top}}}
\newcommand\sbot{_{\text{bot}}}
\newcommand\sll{_{\mathbf l}}
\newcommand\srr{_{\mathbf r}}
\newcommand\geod{\operatorname{\mathbf g}}
\newcommand\mtorus[1]{\boundary U(#1)}
\newcommand\A{\mathbf A}
\newcommand\Aleft[1]{\A\sleft(#1)}
\newcommand\Aright[1]{\A\sright(#1)}
\newcommand\Atop[1]{\A\sbtop(#1)}
\newcommand\Abot[1]{\A\sbot(#1)}
\newcommand\boundvert{{\boundary_{||}}}
\newcommand\storus[1]{U(#1)}
\newcommand\Momega{\omega_M}
\newcommand\nomega{\omega_\nu}
\newcommand\twist{\operatorname{tw}}
\newcommand\modl{M_\nu}
\newcommand\MT{{\mathbb T}}
\newcommand\Teich{{\mathcal T}}
\renewcommand{\Re}{\operatorname{Re}}
\renewcommand{\Im}{\operatorname{Im}}

\title{Quasi-conformal deformations of nonlinearizable germs}

\author{Kingshook Biswas }

\date{}

\thanks{Research partly supported by  Department of Science and Technology research project
grant DyNo. 100/IFD/8347/2008-2009}

\begin{abstract} Let $f(z) = e^{2\pi i \alpha}z + O(z^2), \alpha \in \mathbb{R}$ be a
germ of holomorphic diffeomorphism in $\mathbb{C}$. For $\alpha$
rational and $f$ of infinite order, the space of conformal
conjugacy classes of germs topologically conjugate to $f$ is
parametrized by the Ecalle-Voronin invariants (and in particular
is infinite-dimensional). When $\alpha$ is irrational and $f$ is
nonlinearizable it is not known whether $f$ admits quasi-conformal
deformations. We show that if $f$ has a sequence of repelling
periodic orbits converging to the fixed point then $f$ embeds into
an infinite-dimensional family of quasi-conformally conjugate
germs no two of which are conformally conjugate.

\smallskip

\begin{center}


\end{center}

\end{abstract}

\maketitle

\overfullrule=0pt

\tableofcontents


\section{Introduction}

\medskip

Let $f(z) = e^{2\pi i \alpha}z + O(z^2), \alpha \in
\mathbb{R}/\mathbb{Z}$ be a germ of holomorphic diffeomorphism fixing the
origin in $\mathbb{C}$. We consider the question of when $f$
admits quasi-conformal deformations, i.e. when do there exist
germs $g$ which are quasi-conformally but not conformally
conjugate to $f$? If $f$ is {\it linearizable} (i.e.
analytically conjugate to the rigid rotation $R_{\alpha}(z) =
e^{2\pi i \alpha} z$) then any germ topologically conjugate to $f$ is
linearizable. In the {\it nondegenerate parabolic} case (i.e.
$\alpha = p/q \in \mathbb{Q}, f^q \neq id$), the quasi-conformal
conjugacy class of $f$ contains an infinite dimensional family of
conformal conjugacy classes parametrized by the Ecalle-Voronin invariants
(\cite{ecalle}, \cite{voronin}). In the {\it irrationally indifferent
nonlinearizable} case ($\alpha$ irrational, $f$ not linearizable), it seems
to be unknown whether quasi-conformal deformations are possible.
We show the following:


\medskip

\begin{theorem}\label{def}{Let $f$ be an irrationally indifferent
nonlinearizable germ with a sequence of repelling periodic orbits
accumulating the origin. Then there is a family of quasi-conformal maps
$\{ h_{\Lambda} \} \Lambda \in \mathcal{M}$ parametrized by
sequences $\mathcal{M} = \{ (\lambda_n)_{n \geq 0} : \lambda_n \in \mathbb{C}, |\lambda_n| > 1
\}$ such that all conjugates of $g_{\Lambda} = h_{\Lambda} \circ f \circ
h^{-1}_{\Lambda}$ are holomorphic, and $g_{\Lambda_1}, g_{\Lambda_2}$ are
conformally conjugate if and only if all but finitely
many terms of $\Lambda_1$ and $\Lambda_2$ agree. For fixed $z$,
$h_{\Lambda}(z)$ depends holomorphically on each $\lambda_n$.}
\end{theorem}

\medskip

The proof proceeds as follows: given $\Lambda = (\lambda_n)$ the
germ $g_{\Lambda}$ is obtained by quasi-conformally deforming
$f$ near its periodic orbits. Each periodic orbit is attracting for
$f^{-1}$, and it is possible to construct an $f$-invariant Beltrami differential
on each basin of attraction such that the quasi-conformal map $h_{\Lambda}$
rectifying the Beltrami differential conjugates $f$ to a holomorphic germ
$g_{\Lambda}$ having multipliers $\lambda_n$ at the periodic orbits.
The multipliers at periodic orbits being invariant under
conformal conjugacies, the conclusion of the theorem follows.

\medskip

Examples of germs satisfying the hypothesis of the theorem are germs
of rational maps of degree $d$ with $\alpha$ satisfying the Cremer
condition of degree $d$ (see for example Milnor \cite{milnor}, Ch.
8)
$$
\limsup \frac{\log \log q_{n+1}}{q_n} > \log d
$$
(where $(p_n/q_n)$ are the continued fraction convergents of $\alpha$)
which ensures that the fixed point is accumulated by periodic orbits (only
finitely many of which can be repelling for a rational map).

\medskip

Perez-Marco has shown (\cite{perezmcirclemaps}) for any nonlinearizable
germ $f$ the existence of a unique monotone one-parameter family
$(K_t)_{t > 0}$ of full, totally invariant continua called {\it hedgehogs}
containing the fixed point. In \cite{k1} it is proved that any
conformal mapping in a neighbourhood of a hedgehog $K$ of a germ $f_1$
mapping $K$ to a hedgehog of a germ $f_2$ necessarily conjugates $f_1$ to $f_2$.
As a corollary of Theorem \ref{def} we have

\medskip

\begin{cor}\label{motion}{There exists a holomorphic motion
$\phi : \mathbb{D}^* \times \hat{\mathbb{C}} \to \hat{\mathbb{C}}$
of $\hat{\mathbb{C}}$ over $\mathbb{D}^*$ and a hedgehog $K$
such that all the sets $\phi(t,K)$ are
hedgehogs, all of which are quasi-conformal images of $K$, but for $s \neq
t$, $\phi(s,K)$ cannot be conformally mapped to $\phi(t,K)$ .}
\end{cor}

\medskip

\noindent{\bf Acknowledgements.} The author thanks Ricardo
Perez-Marco for many helpful discussions and comments. This
article was written in part during a visit to the Chennai
Mathematical Institute. The author is very grateful to the organizers for their warmth
and hospitality.

\medskip
%
%


\section{Deformations.} We fix a germ $f(z) = e^{2\pi i \alpha}z + O(z^2), \alpha \in \mathbb{R} - \mathbb{Q}$
 and a neighbourhood $U$ of the origin
such that $f$ and $f^{-1}$ are univalent on a neighbourhood of $\overline{U}$.
By a {\it periodic orbit} or {\it cycle} of $f$ of order $q \geq 1$
we mean a finite set $\mathcal{O} = \{ z_1,\dots,z_q \} \subset U$ such that $f(z_i) =
z_{i+1}, 1 \leq i \leq q-1$ and $f(z_q) = z_1$.
The {\it multiplier} at the periodic orbit is defined to be $\lambda =
f'(z_1)f'(z_2) \dots f'(z_q) = (f^q)'(z_i)$. The periodic orbit is
called {\it attracting}, {\it indifferent} and {\it repelling}
according as $|\lambda| < 1, |\lambda| = 1$ and $|\lambda| > 1$
respectively. A periodic orbit for $f$ of multiplier $\lambda$ is a
periodic orbit for $f^{-1}$ of multiplier $\lambda^{-1}$. The basin of
attraction of an attracting periodic
cycle is defined by $\mathcal{A}(\mathcal{O}, f) := \{ z \in U :
f^n(z) \to \mathcal{O} \hbox{ as } n \to +\infty \}$.

\medskip

We observe that $f$ can have only finitely many cycles of a given order $q$ in $U$
(by the uniqueness principle). Since $f$ is asymptotic to an
irrational rotation, i.e. $f(z)/z \to e^{2\pi i \alpha}$ as $z \to
0$, it follows that if $f$ has small cycles then the orders of the
cycles must go to infinity.

\medskip

\subsection{Deforming repelling periodic orbits.} Given a repelling
periodic orbit $\mathcal{O}$ of $f$ with multiplier $\lambda$ and
$|\lambda'| > 1$ we deform the multiplier of $f$ quasi-conformally
from $\lambda$ to $\lambda'$ by constructing a corresponding $f$-invariant Beltrami differential
$\mu = \mu(\mathcal{O}, \lambda, \lambda')$ on the basin of attraction $\mathcal{A}(\mathcal{O},
f^{-1})$ as follows:

\medskip

By K\"oenigs linearization theorem (see \cite{milnor}, Ch. 6)
for any repelling periodic orbit $\mathcal{O}$ of $f$ with multiplier
$\lambda$ there exists a unique holomorphic map $\phi$ defined on a
neighbourhood of $\mathcal{O}$ such that $\phi(z_i) = 0, \phi'(z_i) =
1, i=1,\dots,n$ and $\phi(f^q(z)) = \lambda
\phi_{\lambda}(z)$. The conformal isomorphism $L : \mathbb{C}^*
\to \mathbb{C}/\mathbb{Z}, w \mapsto \xi = \frac{1}{2\pi i}\log w$
conjugates the linear map $w \mapsto \lambda w$ on $\mathbb{C}^*$
to the translation $\xi \mapsto \xi + \tau$ on
$\mathbb{C}/\mathbb{Z}$ where $\tau = L(\lambda)$ and $\Im \tau < 0$.

\medskip

Let $\tau' = L(\lambda')$ and let $K$ be the real linear map on
$\mathbb{C}$ defined by $K(1) = 1, K(\tau) = \tau'$. Then $K$
commutes with the translation by one and hence gives a quasi-conformal
orientation preserving (since $\Im \tau, \Im \tau' < 0$)
homeomorphism $\tilde{K} : \mathbb{C}/\mathbb{Z} \to
\mathbb{C}/\mathbb{Z}$. The Beltrami differential of $\tilde{K}$
is constant and invariant under translations of
$\mathbb{C}/\mathbb{Z}$, and $\tilde{K}$
conjugates the translation $\xi \mapsto
\xi + \tau$ on $\mathbb{C}/\mathbb{Z}$ to $\xi \mapsto \xi +
\tau'$.

\medskip

We let $\mu$ be the Beltrami differential of $\tilde{K}
\circ L \circ \phi_{\lambda}$ restricted to a small neighbourhood $D$ of $z_1 \in \mathcal{O}$.
The map $\tilde{K} \circ \phi_{\lambda}$ conjugates $f^q$
to the translation $\xi \mapsto \xi + \tau'$. So at points
$z, z'=f^q(z)$ in $D$, $\mu$ satisfies the
invariance condition
$$
\mu(z') \overline{(f^q)'(z)} =
\mu(z) (f^q)'(z)
$$
So $\mu$ extended to the neighbourhood $V = D \cup f(D) \cup \dots \cup
f^{q-1}(D)$ of $\mathcal{O}$ by putting
$$
\mu(f^j(z)) = \mu(z)
\frac{(f^j)'(z)}{\overline{(f^j)'(z)}}
$$
for
$z \in D, j=1,\dots,q-1$ is an $f$-invariant Beltrami
differential. Similarly the above equation allows us to extend
$\mu$ to an $f$-invariant Beltrami differential $\mu(\mathcal{O}, \lambda, \lambda')$ on
$\mathcal{A}(\mathcal{O}, f^{-1}) = \cup_{n \geq 1} f^n(V)$.

\medskip

\begin{lemma}\label{multiplier}{Any quasi-conformal homeomorphism $h$
with Beltrami coefficient equal to $\mu = \mu(\mathcal{O}, \lambda, \lambda')$
on a neighbourhood of $\mathcal{O}$
conjugates $f$ to a map $g = h \circ
f \circ h^{-1}$ with a periodic orbit $h(\mathcal{O})$ of multiplier $\lambda'$. The
dependence of $\mu(\mathcal{O},
\lambda, \lambda')$ on $\lambda'$ is holomorphic.}
\end{lemma}

\medskip

\noindent{Proof:} For $i = 1,\dots, q$ we let $\psi_i$ be the branch
of $\phi^{-1}$ sending $0$ to $z_i$. By construction the map $k$
defined by $k = \psi_i \circ L^{-1} \circ \tilde{K}
\circ L \circ \phi(z)$ for $z$ in a neighbourhood of $z_i$ has Beltrami
coefficient equal to $\mu$ and conjugates
$f$ on a neighbourhood of $\mathcal{O}$ to a holomorphic map
$f_1 = k \circ f \circ k^{-1}$ with periodic orbit $k(\mathcal{O})$
and multiplier $\lambda'$. Since $h$ and $k$ have the same
Beltrami coefficient the map $h \circ k^{-1}$ is holomorphic, and
conjugates $f_1$ to $g$, so the multipliers of $f_1$ and $g$ are
equal. The Beltrami differential of $\tilde{K}$ is constant equal to
$\frac{-\log(\lambda'/\lambda)}{2\log|\lambda| + \log(\lambda'/\lambda)}$
which depends holomorphically on $\lambda'$, so the Beltrami differential
$\mu(\mathcal{O}, \lambda, \lambda')$ of $\tilde{K}
\circ L \circ \phi_{\lambda}$ depends holomorphically on $\lambda'$. $\diamond$

\medskip

\subsection{Deforming germs with small cycles.}

\medskip

We use the Beltrami differentials $\mu(\mathcal{O}, \lambda, \lambda')$
to deform a germ with small cycles as follows:

\medskip

\noindent{Proof of Theorem \ref{def}}: Let $f$ be a germ with a
sequence of repelling small cycles of orders $(q_n)$ (which we may assume
to be strictly increasing).
The basins of attraction
$\mathcal{A}(\mathcal{O}, f^{-1})$ of distinct repelling cycles of
$f$ are disjoint, so for any $\Lambda = (\lambda_n) \in
\mathcal{M}$ we can define an $f$-invariant Beltrami differential
$\mu_{\Lambda}$ on $U$ by putting $\mu_{\Lambda}(z) = \mu(\mathcal{O}, \lambda,
\lambda_n)(z)$ when $z$ belongs to $\mathcal{A}(\mathcal{O}, f^{-1})$ for a
cycle $\mathcal{O}$ of order $q_n$ and multiplier $\lambda$, and
$\mu_{\Lambda}(z) = 0$ otherwise. Let $h_{\Lambda}$ be the unique
quasi-conformal homeomorphism with Beltrami coefficient
$\mu_{\Lambda}$ fixing $0,1,\infty$ given by the Measurable
Riemann Mapping Theorem. Then $g_{\Lambda} = h_{\Lambda} \circ f \circ h^{-1}_{\Lambda}$ is a holomorphic germ
fixing the origin with small cycles. By Naishul's theorem
\cite{naishul} the multiplier at an indifferent fixed point is a topological conjugacy
invariant so $g'_{\Lambda}(0) = e^{2\pi i \alpha}$. By Lemma
\ref{multiplier}, the multiplier of $g_{\Lambda}$ at {\it all} its
repelling cycles of order $q_n$ is equal to $\lambda_n$.

\medskip

If $\phi$ is a holomorphic germ conjugating two such germs $g_{\Lambda_1}, g_{\Lambda_2}$,
$\phi$ must take all repelling cycles of $g_{\Lambda_1}$
of order $q_n$ (in its domain) to repelling cycles of $g_{\Lambda_2}$ of order $q_n$ and
preserve multipliers, so the sequences of multipliers $\Lambda_1, \Lambda_2$
must agree for all but finitely many terms. It follows from Lemma \ref{multiplier}
that $\mu_{\Lambda}$ depends holomorphically on each $\lambda_n$, hence by the
Measurable Riemann Mapping Theorem so does $h_{\Lambda}(z)$ for fixed $z$. $\diamond$

\medskip

\noindent{Proof of Corollary \ref{motion}:} Let $f$ be as above and
$K \subset U$ be a hedgehog of $f$.
For $t \in \mathbb{D}^*$ we let $\Lambda_t$ be the constant
sequence $(\lambda_n = 1/t)$. Then $\mu_{\Lambda_t}$ depends
holomorphically on $t$, and $\phi : (t,z) \mapsto
h_{\Lambda_t}(z)$ gives the required holomorphic motion.
$\diamond$


\bibliography{deform}
\bibliographystyle{alpha}

\medskip

\noindent Ramakrishna Mission Vivekananda University,
Belur Math, WB-711202, India

\noindent email: kingshook@rkmvu.ac.in

\end{document}